\documentclass[11pt,oneside]{amsart}
\title[Infinite time Turing machines and an application to equivalence relations]{Infinite time Turing machines and an application to the
hierarchy of equivalence relations on the reals}
\thanks{The second author's research has been supported in part
by grants from the Research Foundation of CUNY and from the
National Science Foundation.}
\author{Samuel Coskey}
\address{The Fields Institute, 222 College Street, Toronto Ontario M5T
  3J1 \& The York University Department of Mathematics and Statistics,
  N520 Ross, 4700 Keele Street, Toronto, Ontario M3J 1P3}
\email{scoskey@nylogic.org, http://math.rutgers.edu/$\sim$scoskey}
\author{Joel David Hamkins} 
\address{The Graduate Center of The City University of New York,
  Mathematics Program, 365 Fifth Avenue, New York, NY 10016 \& College
  of Staten Island of CUNY, Mathematics, 2800 Victory Boulevard,
  Staten Island, NY 10314}
\email{jhamkins@gc.cuny.edu, http://jdh.hamkins.org}

\usepackage[marginratio=1:1]{geometry}
\usepackage{mathpazo,bm}
\usepackage{setspace}
\usepackage{amssymb}

\newtheorem{theorem}{Theorem}

\newtheorem{question}[theorem]{Question}

\newcommand{\QED}{\end{proof}}

\newcommand{\WO}{{\rm WO}}
\newcommand{\HC}{{\rm HC}}

\newcommand{\R}{{\mathbb R}}
\newcommand{\N}{{\mathbb N}}
\newcommand{\ZFC}{\hbox{$\mathrm ZFC$}}
\newcommand{\intersect}{\cap}
\newcommand{\of}{\subseteq}
\newcommand{\ofnoteq}{\subsetneq}
\newcommand{\proves}{\vdash}
\newcommand{\from}{\mathord{\,\vbox{\baselineskip=2pt\lineskiplimit=0pt
                         \hbox{.}\hbox{.}\hbox{.}}}\;}
\newcommand{\gcode}[1]{\ulcorner\!#1\!\urcorner}
\newcommand{\PP}{\mathrm{P}}
\newcommand{\NP}{\mathrm{NP}}
\newcommand{\PSpace}{\mathrm{PSPACE}}
\newcommand{\coNP}{\mathrm{co}\hbox{-}\mathrm{NP}}
\newcommand{\restrict}{\mathbin{\upharpoonright}}
\newcommand{\elesub}{\prec}
\newcommand{\trianglelt}{\triangleleft}
\newcommand{\set}[1]{\{\,{#1}\,\}}
\newcommand{\singleton}[1]{\left\{{#1}\right\}}
\newcommand{\st}{\mid}

\newcommand{\omegaCK}{{\omega_1^{\hbox{\tiny\sc CK}}}}
\newcommand{\oiso}{\mathord{\cong}}
\newcommand{\oequiv}{\mathord{\equiv}}
\def\<#1>{\langle#1\rangle}

\newcommand{\cell}[1]{\boxit{\hbox to 17pt{\strut\hfil$#1$\hfil}}}
\newcommand{\head}[2]{\lower2pt\vbox{\hbox{\strut\footnotesize\it\hskip3pt#2}\boxit{\cell#1}}}
\newcommand{\boxit}[1]{\setbox4=\hbox{\kern2pt#1\kern2pt}\hbox{\vrule\vbox{\hrule\kern2pt\box4\kern2pt\hrule}\vrule}}
\newcommand{\Col}[3]{\hbox{\vbox{\baselineskip=0pt\parskip=0pt\cell#1\cell#2\cell#3}}}
\newcommand{\tapenames}{\raise 10pt\vbox to 40pt{\hbox to .8in{\it\hfill input: \strut}\vfill\hbox to
.8in{\it\hfill scratch: \strut}\vfill\hbox to .8in{\it\hfill output: \strut}}}
\newcommand{\Head}[4]{\lower2pt\vbox{\hbox to25pt{\strut\footnotesize\it\hfill#4\hfill}\boxit{\Col#1#2#3}}}
\newcommand{\Dots}{\raise 5pt\vbox to 40pt{\hbox{\ $\cdots$\strut}\vfill\hbox{\ $\cdots$\strut}\vfill\hbox{\ $\cdots$\strut}}}
\newcommand{\jump}{{\!\triangledown}}
\newcommand{\Jump}{{\!\blacktriangledown}}
\def\ilt{<_{\infty}}
\def\ileq{\leq_{\infty}}
\def\iequiv{\equiv_{\infty}}
\renewcommand{\th}{{\hbox{\scriptsize th}}}
\newcommand{\df}{\it}
\newcommand{\converges}{\mathord{\downarrow}}

\newcommand{\UnderTilde}[1]{{\setbox1=\hbox{$#1$}\baselineskip=0pt\vtop{\hbox{$#1$}\hbox to\wd1{\hfil$\sim$\hfil}}}{}}
\newcommand{\Undertilde}[1]{{\setbox1=\hbox{$#1$}\baselineskip=0pt\vtop{\hbox{$#1$}\hbox to\wd1{\hfil$\scriptstyle\sim$\hfil}}}{}}
\newcommand{\undertilde}[1]{{\setbox1=\hbox{$#1$}\baselineskip=0pt\vtop{\hbox{$#1$}\hbox to\wd1{\hfil$\scriptscriptstyle\sim$\hfil}}}{}}
\newcommand{\UnderdTilde}[1]{{\setbox1=\hbox{$#1$}\baselineskip=0pt\vtop{\hbox{$#1$}\hbox to\wd1{\hfil$\approx$\hfil}}}{}}
\newcommand{\Underdtilde}[1]{{\setbox1=\hbox{$#1$}\baselineskip=0pt\vtop{\hbox{$#1$}\hbox to\wd1{\hfil\scriptsize$\approx$\hfil}}}{}}

\begin{document}

\begin{abstract}
  We describe the basic theory of infinite time Turing machines and
  some recent developments, including the infinite time degree theory,
  infinite time complexity theory, and infinite time computable model
  theory.  We focus particularly on the application of infinite time
  Turing machines to the analysis of the hierarchy of equivalence
  relations on the reals, in analogy with the theory arising from
  Borel reducibility.  We define a notion of infinite time
  reducibility, which lifts much of the Borel theory into the class
  $\Undertilde\Delta^1_2$ in a satisfying way.
\end{abstract}

\maketitle

\noindent Infinite time Turing machines fruitfully extend the
operation of ordinary Turing machines into transfinite ordinal time
and by doing so provide a robust theory of computability on the
reals. In a mixture of methods and ideas from set theory, descriptive
set theory and computability theory, the approach provides infinitary
concepts of computability and decidability on the reals, which climb
nontrivially into the descriptive set-theoretic hierarchy (at the
level of $\Delta^1_2$) while retaining a strongly computational
nature. With infinite time Turing machines, we have infinitary
analogues of numerous classical concepts, including the infinite time
Turing degrees, infinite time complexity theory, infinite time
computable model theory, and now also the infinite time analogue of
the theory of Borel equivalence relations under Borel reducibility.

In this article, we shall give a brief review of the machines and
their basic theory, and then explain in a bit more detail our recent
application of infinite time computability to an analogue of Borel
equivalence relation theory, a full account of which is given in
\cite{CoskeyHamkins:InfiniteTimeComputableEquivalenceRelations}.  The
basic idea of this application is to replace the concept of Borel
reducibility commonly used in that theory with forms of infinite time
computable reducibility, and study the accompanying hierarchy of
equivalence relations.  This approach retains much of the Borel
analysis and results, while also illuminating a part of the hierarchy
of equivalence relations that seems beyond the reach of the Borel
theory, including many highly canonical equivalence relations that are
infinite time computable but not Borel, such as the isomorphism
relations for diverse classes of countable structures.

Major parts of this article are adapted from the surveys
\cite{Hamkins2007:ASurveyOfInfiniteTimeTuringMachines} and
\cite{Hamkins2005:InfinitaryComputabilityWithITTM} and from our
article
\cite{CoskeyHamkins:InfiniteTimeComputableEquivalenceRelations} on
infinite time computable equivalence relations. Infinite time Turing
machines were first studied by Hamkins and Kidder in 1989, with the
core introduction provided by Hamkins and Lewis
\cite{HamkinsLewis2000:InfiniteTimeTM}.  The theory has now been
extended by many others, including Philip Welch, Peter Koepke,
Benedikt L\"owe, Daniel Seabold, Ralf Schindler, Vinay Deolalikar,
Russell Miller, Steve Warner, Giacomo Lenzi, Erich Monteleone, Samuel
Coskey and others. Numerous precursors to the theory include
Blum-Shub-Smale machines (1980s), B\"uchi machines (1960s) and
accompanying developments, Barry Burd's model of Turing machines with
``blurs" at limits (1970s), the extensive development of
$\alpha$-recursion and $E$-recursion theory, a part of higher
recursion theory (since the 1970s), Jack Copeland's accelerated Turing
machines (1990s), Ryan Bissell-Siders' ordinal machines (1990s), and
more recently, Peter Koepke's ordinal Turing machines and ordinal
register machines (2000s). The expanding literature involving infinite
time Turing machines includes \cite{HamkinsLewis2000:InfiniteTimeTM},
\cite{Welch99:FriedmansTrick}, \cite{Welch2000:Eventually},
\cite{Welch2000:LengthsOfITTM}, \cite{Loewe2001:RevisionSequences},
\cite{HamkinsSeabold2001:OneTape},
\cite{HamkinsLewis2002:PostProblem}, \cite{Schindler:Pnot=NP},
\cite{HamkinsWelch2003:Pf=NPf}, \cite{Hamkins2002:Turing},
\cite{Hamkins2004:SupertaskComputation},
\cite{LenziMonteleone2004:OnFixpointArithmeticAndITTM},
\cite{DeolalikarHamkinsSchindler2005:NPcoNP},
\cite{HamkinsMillerSeaboldWarner2007:InfiniteTimeComputableModelTheory},
\cite{Hamkins2005:InfinitaryComputabilityWithITTM},
\cite{Welch:OnAQuestionOfDeolalikarHamkinsSchindler},
\cite{Welch2005:ActionOfOneTapeMachines},
\cite{Koepke2005:TuringComputationsOnOrdinals},
\cite{Hamkins2007:ASurveyOfInfiniteTimeTuringMachines},
\cite{HamkinsMiller:PostsProblemForORMsExplicitApproach},
\cite{HamkinsMiller2007:PostsProblemForORMs},
\cite{HamkinsLinetskyMiller2007:ComplexityOfQuicklyDecidableORMSets}
and others.

\section{A Brief Review of Infinite time Turing machines}

Infinite time Turing machines have exactly the same hardware as their
classical finite time counterparts, with a head moving back and forth
on a semi-infinite paper tape, writing $0$s and $1$s according to the
rigid instructions of a finite program with finitely many states. What
is new about the infinite time Turing machines is that their operation
is extended into transfinite ordinal time. For convenience, the
machines are implemented with a three-tape model, with separate tapes
for input, scratch work and output.
\begin{figure}[h]
$$\tapenames\Col010\Col010\Head101q\Col110\Col101\Col100\Col011\Col001\Dots$$
\end{figure}
The machine operates at successor stages of computation in exactly the
classical manner, according to the program instructions. Computation
is extended to limit ordinal stages simply by defining the limit
configuration of the machines. The idea is to try to preserve as much
as possible the information that the computation has been creating up
to that stage, preserving it in the limit configuration as a kind of
limit of the earlier configurations. Specifically, at any limit
ordinal stage $\xi$, the machine enters what we call the {\it limit}
state, one of the distinguished states along with the {\it start} and
{\it halt} states; the head is reset to the first cell at the left;
and each cell of the tape is updated with the $\limsup$ of the values
previously displayed in that cell. Having thus specified the complete
configuration of the machine at stage $\xi$, the computation may now
continue to stage $\xi+1$ and so on.  Computational output is given
only when the machine explicitly enters the {\it halt} state, and
computation ceases when this occurs.

Since the tapes naturally accommodate infinite binary strings---and
there is plenty of time for the head to inspect every cell---the
natural context for input and output to the machines is Cantor space
$2^\omega$, which we denote by $\R$ and refer to as the reals. Thus,
the machines provide an infinitary notion of computability on the
reals. A program $p$ computes the partial function
$\varphi_p\from\R\to\R$, defined by $\varphi_p(x)=y$ if program $p$ on
input $x$ yields output $y$, where the output of a computation is the
content of the output tape when the machine enters the {\it halt}
state. A subset $A\of\R$ is \emph{infinite time decidable} if the
characteristic function of $A$ is infinite time computable.  The set
$A$ is \emph{infinite time semi-decidable} if the constant partial
function $1\restrict A$ is computable.  This is equivalent to $A$
being the domain of an infinite time computable function (but not
necessarily equivalent to $A$ being the range of such a
function). Elementary results in
\cite{HamkinsLewis2000:InfiniteTimeTM} show that the arithmetic sets
are exactly those that are decidable in time uniformly less than
$\omega^2$ and the hyperarithmetic sets are those that are decidable
in time less than some recursive ordinal. The power of the machines,
however, reaches much higher than this into the descriptive set
theoretic hierarchy.

For example, every $\Pi^1_1$ and $\Sigma^1_1$ set is infinite time
decidable. To see this, it suffices to show that the complete
$\Pi^1_1$ set $\WO$, consisting of reals coding a well-ordered
relation on $\omega$, is infinite time computable. This is
accomplished by the {\it count-through} argument of \cite[Theorem
2.2]{HamkinsLewis2000:InfiniteTimeTM}, which we should like to sketch
here. Given a real $x$, we view it as coding the relation
$\trianglelt$ on $\omega$ for which $n\trianglelt m$ if and only if
the $\<n,m>$ bit of $x$ is $1$. The assertion that $\trianglelt$ is a
linear order is arithmetic in $x$, and therefore easily determined by
the machines. After this, the machine will check for well-foundedness
essentially by counting through the order, relying on the fact that
the computational steps are themselves well-ordered. Specifically, the
machine places an initial guess for the current minimal element in the
relation $\trianglelt$, updating it with better guesses as they are
encountered. At each revision, the machine flashes a certain master
flag, so that at the limit stage the machine can know if the guess was
changed infinitely often, indicating ill-foundedness (the machine
should reset the master flag at limits of limit stages). Otherwise,
the true current minimal element has been found, and so the machine
can delete all mention of it from the field of the relation coded by
$x$. Iterating this, the algorithm in effect systematically erases the
well-founded initial segment of the relation coded by the input real,
until either nothing is left or the ill-founded part is discovered,
either of which can be determined. In this way, membership in $\WO$ is
infinite time decidable. It follows that every $\Pi^1_1$ and
$\Sigma^1_1$ set is infinite time decidable, and so the machines climb
properly into $\Delta^1_2$. Meanwhile, the class of infinite time
decidable sets is easily observed to be contained in $\Delta^1_2$, and
in fact the class $\Delta^1_2$ is closed under the infinite time jump
operations and is therefore stratified by a significant part of the
infinite time Turing degrees.

Although transfinite, computations are nevertheless inherently
countable, since an easy cofinality argument establishes that every
computation either halts or repeats by some countable ordinal
stage. An ordinal $\alpha$ is said to be {\df clockable}, if there is
a computation $\varphi_p(0)$ halting on exactly the $\alpha^\th$ step.
A real $x$ is {\df writable} if it is the output of a computation
$\varphi_p(0)$, and an ordinal is writable if it is coded by such a
real. Because there are only countably many programs, it follows that
there are only countably many clockable and writable ordinals. The
clockable and writable ordinals extend through all the recursive
ordinals and far beyond; their supremum is recursively inaccessible
and more. The writable ordinals form an initial segment of the
ordinals, since whenever an ordinal is writable, then the algorithm
writing it can be easily modified to write a code for any smaller
ordinal.  But the same is not true for the clockable ordinals; in the
midst of the clockable ordinals, there are increasingly complex
forbidden regions at which no (parameter-free) infinite time Turing
machine can halt.

Let us quickly sketch the argument that such gaps in the clockable
ordinals exist, since this is an interesting exercise in ordinal
reflection that constitutes a basic method of many later constructions
in the theory. Consider the algorithm that simulates all programs on
input $0$ simultaneously, by some bookkeeping method that reserves and
manages sufficient separate space for each, simulating $\omega$ many
steps of computation for each program in each $\omega$ many steps of
actual computation. Our algorithm might keep careful track of which
programs have halted, and pay attention to find a stage at which none
of the programs halt. Since such a stage exists above the supremum of
all clockable ordinals, we will definitely find such a stage
eventually. Since our algorithm can recognize the first such stage, we
can arrange that it halts immediately after this discovery. So we have
described a computational procedure that will halt at an ordinal stage
that is larger than a stage at which no computations halted, and so
there are gaps in the clockable ordinals, as desired. A careful
analysis of the algorithm shows that the first gap after any clockable
ordinal has order type $\omega$, essentially because it takes $\omega$
many additional steps to realize that a gap has been reached. Modified
algorithms search for longer gaps and show that there must be
increasingly complex gaps at increasingly complex admissible limit
stages---for any clockable or writable ordinal $\alpha$, there are
gaps of size at least $\alpha$. The structure of these gaps exhibits
the same complexity as the infinite time halting problem.

Although it was established in \cite{HamkinsLewis2000:InfiniteTimeTM}
that the clockable and writable ordinals have the same order type,
perhaps the main question left open in that paper was whether the
supremum of these ordinals was the same. This was settled in the
affirmative by Philip Welch in \cite{Welch2000:LengthsOfITTM}. Another
way to describe the result is that whenever program $p$ on input $x$
yields a halting computation, then there is another computation that
writes out a certificate of this computation, a real coding the entire
computation history including a well-ordered relation whose order type
is the length of the computation.  This important fact, far from
obvious, relies on a subtle treatment of eventual writability and
constitutes a foundation of many further developments of the theory,
including the applications we mention in this article.

The reflective aspect of the count-through argument described above
consists of the observation that any decidable property that holds of
a real that might be encountered during the course of a computation
must hold of a writable real, since we may embark on the computational
search to find such a witness and output it when it is found. This
idea is greatly extended by the $\lambda$-$\zeta$-$\Sigma$ theorem of
Philip Welch.  Specifically, \cite{HamkinsLewis2000:InfiniteTimeTM}
defines that a real $x$ is {\df eventually writable} if there is a
computation $\varphi_p(0)$ for which $x$ appears on the output tape
from some point on (even if the computation does not halt), and $x$ is
{\df accidentally writable} if it appears on any of the tapes at any
stage during a computation $\varphi_p(0)$. By coding ordinals with
reals, we obtain the notions of eventually and accidentally writable
ordinals. If $\lambda$ is the supremum of the clockable or writable
ordinals, $\zeta$ is the supremum of the eventually writable ordinals
and $\Sigma$ is the supremum of the accidentally writable ordinals,
then \cite{HamkinsLewis2000:InfiniteTimeTM} establishes
$\lambda<\zeta<\Sigma$. The $\lambda$-$\zeta$-$\Sigma$ theorem of
Welch \cite{Welch2000:Eventually} asserts moreover that
$L_\lambda\elesub_{\Sigma_1} L_\zeta\elesub_{\Sigma_2} L_\Sigma$,
using the initial segments of G\"odel's constructible universe, and
furthermore, that these ordinals are characterized as the least
example of this pattern. This result precisely expresses the sense in
which the algorithms may pull down witnesses from the accidentally
writable realm into the eventually writable or writable realms. At the
heart of the proof and the result is the fact that every computation
repeats the stage $\zeta$ configuration at stage $\Sigma$.

Many of the fundamental constructions of classical finite time
computability theory carry over to the infinite time context. For
example, one can prove the infinite time analogues of the
$smn$-theorem, the Recursion theorem and the undecidability of the
infinite time halting problem, by essentially the classical
arguments. Some other classical facts, however, do not directly
generalize. For example, it is not true in the infinite time context
that if the graph of a function $f$ is semi-decidable, then the
function is computable. This is a consequence of the following:

\begin{theorem}[Lost Melody Theorem]
  There is a real $c$ such that $\singleton{c}$ is infinite time
  decidable, but $c$ is not writable.
\end{theorem}

The real $c$, a lost melody that you cannot sing on your own, although
you can recognize it yes-or-no when someone sings it to you, exhibits
sufficient internal structure that $\singleton{c}$ is decidable, but
is too complicated itself to be writable. That is, we can recognize
whether a given real $y$ is $c$ or not, but we cannot produce $c$ from
nothing. The function $f(x)=c$ with constant value $c$, therefore, is
not computable, because $c$ is not writable, but the graph is
decidable, because we can recognize whether a pair has the form
$(x,c)$.

The infinite time analogue of the halting problem breaks into
lightface and boldface versions, $h=\set{p\st\varphi_p(p)\converges}$
and $H=\set{(p,x)\st\varphi_p(x)\converges}$, respectively.  These are
both semi-decidable and not decidable, but in the infinitary context,
they are not computably equivalent.

The notion of oracle computation lifts to the infinitary context and
gives rise to a theory of relative computability and a rich structure
of degrees. In contrast to the classical theory on $\N$, however, in
the infinite time context we have two natural sorts of oracles to be
used in oracle computations, corresponding to the second order nature
of the theory. First, one can use an individual real as an oracle in
exactly the classical manner, by adjoining an oracle tape on which the
values of that real are written out. This amounts to fixing a
supplemental input parameter and can be viewed as giving rise to a
boldface theory of infinitary computability, just as one allows
arbitrary real parameters in the descriptive set-theoretic treatment
of boldface $\Undertilde\Delta^1_1$ and $\Undertilde\Pi^1_1$. (We
shall explicitly adopt such a boldface perspective in our application
to the theory of equivalence relations under infinite time
reducibility.)  Second, however, one naturally wants somehow to use a
{\it set} of reals as an oracle, although we cannot expect in general
to write such a set out on the tape (perhaps it is even
uncountable). Instead, the oracle tape is empty at the start of
computation, and during the computation the machine may freely write
on this tape; whenever the algorithm calls for it, the machine may
make a membership query about whether the real currently written on
the oracle tape is a member of the oracle or not. Thus, the machine is
able to know of any real that it can produce, whether the real is in
the oracle set or not.

Such oracle computations give rise to a notion of relative
computability $\varphi_p^A(x)$ and therefore a notion of infinite time
omputable reduction $A\ileq B$ and the accompanying infinite time
degree relation $A\iequiv B$.  For any set $A$, we have the lightface
jump $A^\jump$ and the boldface jump $A^\Jump$, corresponding to the
two halting problems, relativized to $A$. The boldface jump jumps much
higher than the lightface jump, as
\cite{HamkinsLewis2000:InfiniteTimeTM} establishes that $A\ilt
A^\jump\ilt A^\Jump$, as well as $A^{\jump\Jump}\iequiv A^\Jump$ and a
great number of other interesting interactions. The infinite time
analogue of Post's problem, the question of whether there are
intermediate semi-decidable degrees between $0$ and the jump
$0^\jump$, was settled by \cite{HamkinsLewis2002:PostProblem} in an
answer that cuts both ways:

\begin{theorem} The infinite time analogue of Post's problem has both
  affirmative and negative solutions.
  \begin{enumerate}
  \item There are no reals $z$ with $0\ilt z\ilt 0^\jump$.
  \item There are sets of reals $A$ with $0\ilt A\ilt
    0^\jump$. Indeed, there are incomparable semi-decidable sets of
    reals $A\perp_\infty B$.
  \end{enumerate}
\end{theorem}

The degrees of the accidentally writable reals are linearly ordered
and in fact form a well-ordered hierarchy of order type $\zeta+1$,
which corresponds also to their order of earliest appearance on any
computation. In other work, Welch \cite{Welch99:FriedmansTrick} found
minimality in the infinite time Turing degrees. Hamkins and Seabold
\cite{HamkinsSeabold2001:OneTape} analyzed one-tape versus multi-tape
infinite time Turing machines, and Benedikt L\"owe
\cite{Loewe2001:RevisionSequences} observed the connection between
infinite time Turing machines and revision theories of truth.

\section{Some applications and extensions}

Let us briefly describe a few of the recent developments and
extensions of infinite time Turing machines, such as the rise of
infinite time complexity theory and the introduction of infinite time
computable model theory.  After this, in the following section we
shall go into greater detail concerning the application of infinite
time Turing machines to an analogue of the theory of Borel equivalence
relations.

Ralf Schindler \cite{Schindler:Pnot=NP} initiated the study of
infinite time complexity theory by solving the infinite time Turing
machine analogue of the $\PP$ versus $\NP$ question. To define the
polynomial class $\PP$ in the infinite time context, Schindler
observed simply that all reals have length $\omega$ and the polynomial
functions of $\omega$ are bounded by those of the form
$\omega^n$. Thus, he defined that a set $A\of\R$ is in $\PP$ if there
is a program $p$ and a natural number $n$ such that $p$ decides $A$
and halts on all inputs in time before $\omega^n$. The set $A$ is in
$\NP$ if there is a program $p$ and a natural number $n$ such that
$x\in A$ if and only if there is $y$ such that $p$ accepts $(x,y)$,
and $p$ halts on all inputs in time less than $\omega^n$. Schindler
proved $\PP\neq\NP$ for infinite time Turing machines in
\cite{Schindler:Pnot=NP}, using methods from descriptive set theory to
analyze the complexity of the classes $\PP$ and $\NP$. This has now
been generalized in joint work
\cite{DeolalikarHamkinsSchindler2005:NPcoNP} to the following, where
the class $\coNP$ consists of the complements of sets in $\NP$.

\begin{theorem}
  $\PP\neq\NP\intersect\coNP$ for infinite time Turing machines.
\end{theorem}

This proof appears in \cite{DeolalikarHamkinsSchindler2005:NPcoNP}. It
follows that $\PP\neq\NP$ for infinite time Turing machines. (This
result has no bearing whatsoever on the finitary classical $\PP\neq
\NP$ question.) Some of the structural reasons behind
$\PP\neq\NP\intersect\coNP$ are revealed by placing the classes $\PP$
and $\NP$ within a larger hierarchy of complexity classes $\PP_\alpha$
and $\NP_\alpha$ using computations of size bounded below
$\alpha$. Results in \cite{DeolalikarHamkinsSchindler2005:NPcoNP}
showed, for example, that the classes $\NP_\alpha$ are identical for
$\omega+2\leq\alpha\leq\omegaCK$, but nevertheless,
$\PP_{\alpha+1}\ofnoteq P_{\alpha+2}$ for any clockable limit ordinal
$\alpha$. It follows, since the $\PP_\alpha$ are steadily increasing
while the classes $\NP_\alpha\intersect\coNP_\alpha$ remain the same,
that $\PP_\alpha\ofnoteq\NP_\alpha\intersect\coNP_\alpha$ for any
ordinal $\alpha$ with $\omega+2\leq\alpha<\omegaCK$.  Thus,
$\PP\neq\NP\intersect\coNP$. Nevertheless, we attain equality at the
supremum $\omegaCK$ with
$$\PP_{\omegaCK}=\NP_{\omegaCK}\intersect\coNP_{\omegaCK}.$$
In fact, this is an instance of the equality
$\Delta^1_1=\Sigma^1_1\intersect\Pi^1_1$, and one can thereby begin to
see how the theory of infinite time Turing machines grows naturally
into descriptive set theory.

This same pattern of inequality
$\PP_\alpha\ofnoteq\NP_\alpha\intersect\coNP_\alpha$ is mirrored
higher in the hierarchy, whenever $\alpha$ lies strictly within a
contiguous block of clockable ordinals, with the corresponding
$\PP_\beta=\NP_\beta\intersect\coNP_\beta$ for any $\beta$ that begins
a gap in the clockable ordinals. In addition, the question is settled
in \cite{DeolalikarHamkinsSchindler2005:NPcoNP} for the other
complexity classes $\PP^+$, $\PP^{++}$ and $\PP^f$.  Benedikt L\"owe
has introduced analogues of $\PSpace$.

The subject of infinite time computable model theory was introduced in
\cite{HamkinsMillerSeaboldWarner2007:InfiniteTimeComputableModelTheory}.
Computable model theory is model theory with a view to the
computability of the structures and theories that arise.  Infinite
time computable model theory carries out this program with the notion
of infinite time computability provided by infinite time Turing
machines. The classical theory began decades ago with such topics as
computable completeness (Does every decidable theory have a decidable
model?) and computable categoricity (Does every isomorphic pair of
computable models have a computable isomorphism?), and the field has
now matured into a sophisticated analysis of the complexity spectrum
of countable models and theories.

The motivation for a broader context is that, while classical
computable model theory is necessarily limited to countable models and
theories, the infinitary computability context allows for uncountable
models and theories, built on the reals. Many of the computational
constructions in computable model theory generalize from structures
built on $\mathbb N$, using finite time computability, to structures
built on $\R$, using infinite time computability. The uncountable
context opens up new questions, such as the infinitary computable
L\"owenheim-Skolem Theorem, which have no finite time
analogue. Several of the most natural questions turn out to be
independent of \ZFC.

In joint work
\cite{HamkinsMillerSeaboldWarner2007:InfiniteTimeComputableModelTheory},
we defined that a model ${\mathcal A}=\<A,\ldots>$ is infinite time
{\df computable} if $A\of\R$ is decidable and all functions, relations
and constants are uniformly infinite time computable from their
G\"odel codes and input. The structure $\mathcal A$ is {\df decidable}
if one can compute whether ${\mathcal A}\models\varphi[\bar a]$ given
$\gcode{\varphi}$ and $\bar a$. A theory $T$ is infinite time {\df
  decidable} if the relation $T\proves\varphi$ is computable in
$\gcode{\varphi}$.  Because we want to treat uncountable languages,
the natural context for G\"odel codes is $\R$ rather than $\mathbb N$.

The initial question, of course, is the infinite time computable
analogue of the Completeness Theorem: Does every consistent decidable
theory have a decidable model? The answer turns out to be independent
of \ZFC.

\begin{theorem}[\cite{HamkinsMillerSeaboldWarner2007:InfiniteTimeComputableModelTheory}]
  The infinite time computable analogue of the Completeness Theorem is
  independent of \ZFC.
  Specifically:\label{Theorem.ITCCTisIndependent}
  \begin{enumerate}
  \item If\/ $V=L$, then every consistent infinite time decidable
    theory has an infinite time decidable model, in a computable
    translation of the language.
  \item It is relatively consistent with \ZFC\ that there is an
    infinite time decidable theory, in a computably presented
    language, having no infinite time computable or decidable model in
    any translation of the language.
  \end{enumerate}
\end{theorem}

The proof of (1) uses the concept of a {\df well-presented} language
$\mathcal L$, for which there is an enumeration of the symbols
$\<s_\alpha\st\alpha<\delta>$ such that from any $\gcode{s_\alpha}$
one can uniformly compute a code for the prior symbols
$\<\gcode{s_\beta}\st\beta\leq\alpha>$.  One can show that every
consistent decidable theory in a well-presented language has a
decidable model, and if $V=L$, then every computable language has a
well presented computable translation. For (2), one uses the theory
$T$ extending the atomic diagram of $\<\WO,\equiv>$ while asserting
that $f$ is a choice function on the $\equiv$ classes. This is a
decidable theory, but for any computable model ${\mathcal
  A}=\<A,\equiv,f>$ of $T$, the set $\set{f(c_u)\st u\in\WO}$ is
$\Sigma^1_2$ and has cardinality $\omega_1$. It is known to be
consistent with \ZFC\ that no $\Sigma^1_2$ set has size $\omega_1$.

For the infinite time analogues of the L\"owenheim-Skolem Theorem, we
proved for the upward version that every well presented infinite time
decidable model has a proper elementary extension with a decidable
presentation, and for the downward version, every well presented
uncountable decidable model has a countable decidable elementary
substructure. There are strong counterexamples to a full direct
generalization of the L\"owenheim-Skolem theorem, however, because
\cite{HamkinsMillerSeaboldWarner2007:InfiniteTimeComputableModelTheory}
provides a computable structure $\<\R,U>$ on the entire set of reals,
which has no proper computable elementary substructure.

Some of the most interesting work involves computable quotients. A
structure has an infinite time computable {\df presentation} if it is
isomorphic to a computable structure, and has a computable {\df
  quotient presentation} if it is isomorphic to the quotient of a
computable structure by a computable equivalence relation (a
congruence). For structures on $\mathbb N$, in either the finite or
infinite time context, these notions are equivalent, because one can
computably find the least element of any equivalence class. For
structures on $\R$, however, computing such distinguished elements of
every equivalence class is not always possible.

\begin{question}
  Does every structure with an infinite time computable quotient
  presentation have an infinite time computable
  presentation?\label{Question.QuotientPresentation}
\end{question}

In the finite time theory, or for structures on $\mathbb N$, the
answer of course is Yes. But in the full infinite time context for
structures on $\R$, the answer depends on the set theoretic
background.

\begin{theorem}\label{Theorem.QuotientPresentation}
  The answer to Question \ref{Question.QuotientPresentation}
  is independent of \ZFC. Specifically, %
  \begin{enumerate}
  \item It is relatively consistent with \ZFC\ that every structure
    with an infinite time computable quotient presentation has an
    infinite time computable presentation.
  \item It is relatively consistent with \ZFC\ that there is a
    structure having an infinite time computable quotient
    presentation, but no infinite time computable presentation.
  \end{enumerate}
\end{theorem}

Let us briefly sketch some of the ideas appearing in the proof. In
order to construct an infinite time computable presentation of a
structure, given a computable quotient presentation, we'd like somehow
to select a representative from each equivalence class, in a
computably effective manner, and build a structure on these
representatives.  Under the set theoretic assumption $V=L$, we can
attach to the $L$-least member of each equivalence class an escort
real that is powerful enough to reveal that it is the $L$-least member
of its class, and build a computable presentation out of these
escorted pairs of reals. (In particular, the new presentation is not
built out of mere representatives from the original class, since these
reals may be too weak; they need the help of their escorts.)  Thus, if
$V=L$, then every structure with a computable quotient presentation
has a computable presentation. On the other side of the independence,
we prove statement 2 by the method of forcing. The structure
$\<\omega_1,<>$ always has a computable quotient presentation built
from reals coding well orders, but there are forcing extensions in
which no infinite time computable set has size $\omega_1$, on
descriptive set theoretic grounds. In these extensions, therefore,
$\<\omega_1,<>$ has a computable quotient presentation, but no
computable presentation.

Let us also briefly discuss some of the alternative models of ordinal
computation to which infinite time Turing machines have given
rise. Peter Koepke \cite{Koepke2005:TuringComputationsOnOrdinals}
introduced the {\df Ordinal Turing Machines}, which generalize the
infinite time Turing machines by extending the tape to transfinite
ordinal length. The limit rules are accordingly adjusted so that the
machine can make use of this extra space. Specifically, rather than
using a special {\it limit} state, the ordinal Turing machines simply
have a fixed order on their (finitely many) states, and at any limit
stage, the state is defined to be the $\liminf$ of the prior
states. The head position is then defined to be the $\liminf$ of the
head positions when the machine was previously in that resulting limit
state. For uniformity, then, Koepke defines that the cells of the tape
use the $\liminf$ of the prior cell values (rather than $\limsup$ as
with the infinite time Turing machines). If the head moves left from a
cell at a limit position, then it appears all the way to the left on
the first cell.

These machines therefore provide a model of computation for functions
on the ordinals, and notions of decidability for classes of
ordinals. The main theorem is that the power of these machines is
essentially the same as that of G\"odel's constructible universe.

\begin{theorem}[Koepke]
  The sets of ordinals that are ordinal Turing machine decidable, with
  finitely many ordinal parameters, are exactly the sets of ordinals
  in G\"odel's constructible universe $L$.
\end{theorem}

Several other infinitary models of ordinal computation are based on a
concept of ordinal registers, and have given rise to a rich
theory. See \cite{Koepke2005:TuringComputationsOnOrdinals},
\cite{Koepke-Siders}, \cite{koepke-ordinal}, \cite{koepke-recursion},
\cite{koepke-itrm}, \cite{HamkinsMiller2007:PostsProblemForORMs},
\cite{HamkinsMiller:PostsProblemForORMsExplicitApproach}, and
\cite{HamkinsLinetskyMiller2007:ComplexityOfQuicklyDecidableORMSets}.

\section{Infinite time computable equivalence relation
theory}

Recently, we have introduced the natural analogue of Borel equivalence
relation theory in which infinite time decidable relations are
compared with respect to infinite time computable reduction functions.
This is motivated in part by the occasional need in the study of Borel
equivalence relations to go beyond Borel.  Indeed, a more powerful
notion of reducibility may be able to accurately compare more complex
relations.  In particular, we shall be able to consider the new
relations which arise out of the infinite time complexity classes.

We begin with a quick introduction to the study of Borel equivalence
relations.  The name of the subject is somewhat of a misnomer---in
fact the principle objects of study are arbitrary equivalence
relations on \emph{standard Borel spaces}, that is, sets equipped with
the Borel structure of a complete separable metric space.  In
applications, we think of an equivalence relation as representing a
classification problem from some other area of mathematics.  For
instance, since any group with domain $\N$ is determined by its
multiplication function, studying the classification problem for
countable groups amounts to studying the isomorphism equivalence
relation on a suitable subspace of $2^{\N\times\N\times\N}$.  For many
more examples, see Section~1.2 of \cite{thomas_workshop}.

The theory of Borel equivalence relations revolves around the
following key notion of complexity.  If $E,F$ are equivalence
relations on standard Borel spaces $X,Y$, then following
\cite{friedman} and \cite{hjorthkechris} we say that $E$ is
\emph{Borel reducible} to $F$, written $E\leq_BF$, iff there exists a
Borel function $f\colon X\rightarrow Y$ such that
\begin{equation}
  \label{reduction}
  x\mathrel{E}x'\iff f(x)\mathrel{F}f(x')\;.
\end{equation}
Borel reducibility measures the complexity of equivalence relations
not as sets of pairs, but as \emph{classification problems}.  That is,
if $E$ is Borel reducible to $F$, then the classification of elements
of $X$ up to $E$ is no harder than the classification of elements of
$Y$ up to $F$.  The by now classical and highly successful study of
Borel equivalence relations consists in part of two major endeavors.
First, one wishes to map out the relationships between numerous
well-understood and naturally occurring equivalence relations.
Second, given a real-life classification problem one should measure
its complexity by comparing it against the mapped-out benchmark
relations.

Some definability condition on the reduction functions (in this case
that they be Borel) is necessary.  Indeed, without any such
restriction reducibility would always be determined by cardinalities
alone.  However, there are cases of natural classifications by
invariants which cannot be computed by a Borel reduction function.
For instance, it is $\Undertilde\Delta^1_2$ and not Borel to compute
the classical Ulm invariants for a countable torsion abelian group.
One might be tempted to form a theory of $\Undertilde\Delta^1_2$
reducibility, but it turns out this notion is too generous.  Indeed,
as we shall see below in Theorem~\ref{Theorem.SemiCollapse}, it may
lump most equivalence relations together into one trivial complexity
class.

We will consider here reduction functions which are computable by an
infinite time Turing machine (see
\cite{CoskeyHamkins:InfiniteTimeComputableEquivalenceRelations} for a
more complete exposition).  Thus, for any two equivalence relations
$E,F$ on $\R$, we say that $E$ is \emph{infinite time computably
  reducible} to $F$, written $E\leq_c F$, if there is an infinite time
computable function $f$ (freely allowing real parameters) satisfying
Equation~\eqref{reduction}.  Similarly, we say that $E$ is
\emph{eventually reducible} to $F$, written $E\leq_eF$, if there is an
eventually computable function $f$ satisfying
Equation~\eqref{reduction}.  Note here that since all uncountable
standard Borel spaces are Borel isomorphic, we lose no generality by
restricting ourselves to equivalence relations with domain $\R$.

Of course, by the remarks in Section~1 (and again emphasizing that we
have allowed parameters) the infinite time computable reductions
include all of the Borel reductions.  Thus, our theory will extend the
classical theory.  Conversely, many classical proofs of
\emph{non}-reducibility $E\not\leq_BF$ rely on methods such as
measure, category, or forcing.  Hence, they frequently ``overshoot''
and show that there does not exist a reduction from $E$ to $F$ which
is Lebesgue measurable, Baire measurable, or absolutely
$\Undertilde\Delta^1_2$ (discussed below), respectively.  Since the
infinite time computable and eventually computable functions enjoy all
three of these properties, it follows in each of these cases that
$E\not\leq_cF$ and even $E\not\leq_eF$, and hence not too much is
``collapsed'' when we pass from the $\leq_B$ hierarchy to the $\leq_c$
and $\leq_e$ hierarchies.

The infinite time notions of reducibility are very closely related to
that of absolutely $\Undertilde\Delta^1_2$ reducibility, which has
been treated in the literature by Hjorth and others.  Recall that a
subset $A\of\R$ is said to be \emph{absolutely}
$\Undertilde\Delta^1_2$ if it is defined by equivalent
$\Undertilde\Sigma^1_2$ and $\Undertilde\Pi^1_2$ formulas which remain
equivalent in every forcing extension.  A function
$f\colon\R\rightarrow\R$ is said to be absolutely
$\Undertilde\Delta^1_2$ if its diagram $\set{(x,n)\st f(x)\in B_n}$ is
absolutely $\Undertilde\Delta^1_2$ (here, $B_n$ runs through the basic
open subsets of $\R$).  We know of very few naturally occurring cases
in which there is an absolutely $\Undertilde\Delta^1_2$ reduction
between two equivalence relations but not an infinite time computable
reduction.  And when there is an infinite time computable reduction,
one can demonstrate that this is the case by simply ``coding up'' an
algorithm which implements the witnessing reduction function.  This
computational approach may be more satisfying than abstractly defining
a reduction function and verifying that it is $\Undertilde\Delta^1_2$
in all forcing extensions.  On the other hand, we do not have any
general tools for establishing non-reducibility by infinite time
computable functions beyond the already established tools mentioned
above, all of which establish non-reducibility by absolutely
$\Undertilde\Delta^1_2$ functions already.  A brief summary of results
due to Hjorth and Kanovei which establish non-reducibility for
absolutely $\Undertilde\Delta^1_2$ functions can be found in Section~5
of \cite{CoskeyHamkins:InfiniteTimeComputableEquivalenceRelations}.
Some deeper results on this notion of reducibility can be found Hjorth
in Chapter~9 of \cite{hjorthbook}.

For an example of ``coding up'' a new (non-Borel) reduction function,
consider the $E_{ck}$ relation defined by $x\mathrel{E}_{ck}y$ if $x$
and $y$ compute (in the ordinary sense) the same ordinals.  We will
compare it against the relation $\oiso_\WO$, which is just the
isomorphism relation restricted to the set of codes for well-orders.
These two relations are not comparable by Borel reductions;
nevertheless they are closely related and this is made precise by the
following result.

\begin{theorem}
  $E_{ck}$ and $\oiso_\WO$ are infinite time computably bireducible.
\end{theorem}

For instance, there is an intuitive reduction from $E_{ck}$ to
$\oiso_\WO$---namely, map $x$ to a code for the supremum of the
ordinals which are computable (in the ordinary sense) from $x$.  And
indeed, this intuition easily translates into a program for an
infinite time Turing machine.  Briefly, the program simply simulates
all ordinary Turing computations, and inspects the real enumerated by
each.  Whenever one of these reals is seen to be code for a
well-order, this code is added to a list.  Finally, the program
computes and outputs a code for the supremum of the ordinals in its
list.

Another obvious benefit to using infinite time computable and
eventually computable reductions is that they are tailor-made to
handle equivalence relations which arise in the study of infinite time
complexity classes.  As a very simple example, consider two of the
most important such equivalence relations: the infinite time degree
relation $\oequiv_\infty$ which was introduced in Section~1, and the
(light face) jump equivalence relation defined by $x\mathrel{J}y$ if
and only if $x^\jump\iequiv y^\jump$.  We have the following (somewhat
trivial) relationship between the two.

\begin{theorem}
  $J$ is eventually reducible to $\oequiv_\infty$ by the function
  which computes the infinite time jump of a real.
\end{theorem}

The program which witnesses this simply simulates all infinite time
programs on input $x$, and whenever one of them halts adds its index
to a list on its output tape.  Since all programs which will halt do
so by stage $\lambda$, the output tape will eventually show $x^\jump$.

Meanwhile, the next result gives a sampling of non-reducibility
results which can be obtained using the methods of Hjorth and Kanovei
discussed above.  Here $=$ of course denotes the equality relation on
$\R$, and $E_0$ the almost equality relation defined by
$x\mathrel{E}_0y$ if and only if $x(n)=y(n)$ for almost all $n$.
Next, $\oiso_\HC$ denotes the isomorphism relation restricted to the
set of codes for hereditarily countable sets.  Finally,
$E_\mathrm{set}$ denotes the relation defined by
$x\mathrel{E}_\mathrm{set}y$ if $x$ and $y$, thought of as codes for
countable sequences of reals, enumerate the same set.

\begin{theorem}
  \label{Theorem.Non-reductions}\
  \begin{enumerate}
  \item $E_0$ does not infinite time computably reduce to $=$.
  \item $E_\mathrm{set}$ does not infinite time computably reduce to
    $E_0$.
  \item $\oiso_\HC$ and $E_\mathrm{set}$ do not infinite time
    computably reduce to $\oiso_\WO$.
  \end{enumerate}
\end{theorem}

Without strong set-theoretic hypotheses, such results cannot be
obtained for reduction functions which are much more general than the
absolutely $\Undertilde\Delta^1_2$ functions.  For instance, the
infinite time semi-computable reduction functions are still well
inside the class $\Undertilde\Delta^1_2$, but if we were to allow
reduction functions in this class, then all of the equivalence
relations in Theorem~\ref{Theorem.Non-reductions} would be reducible
to the equality relation.

\begin{theorem}
  \label{Theorem.SemiCollapse}
  If $V=L$, then every infinite time computable equivalence relation
  on $\R$ is reducible to the equality relation by an infinite time
  semi-computable function.
\end{theorem}

The proof of Theorem~\ref{Theorem.SemiCollapse} uses the same ideas as
in the proof of Theorem~\ref{Theorem.QuotientPresentation}, and as in
that argument, the reduction functions are not selectors for the
relation.  On the other hand, under suitable determinacy hypotheses,
every infinite time semi-computable function is Lebesgue measurable.
In this situation, infinite time semi-computable reducibility again
resembles the more concrete reducibility notions.

We have seen that by expanding the class of reduction functions
available, we are sometimes able to bring a wider class of equivalence
relations under consideration.  A major example of this is the
following generalization of the class of countable Borel equivalence
relations.  Here, a Borel equivalence relation is said to be
\emph{countable} iff every equivalence class is countable.  The
countable relations have become one of the most important collections
studied in the classical theory, since many natural relations lie at
this level and some basic progress has been made in uncovering their
structure under $\leq_B$.  For instance, by a classical result of
Silver, the equality relation $=$ is the $\leq_B$-least countable
Borel equivalence relation.  Moreover, by a deep result of
Kechris-Harrington-Louveau, $E_0$ is the $\leq_B$-least Borel
equivalence relation which is not reducible to $=$.  Thirdly, we have
that there is a $\leq_B$-greatest countable Borel equivalence
relation, denoted $E_\infty$.  The remaining countable Borel
equivalence relations lie in the interval $(E_0,E_\infty)$, and a
result of Adams-Kechris implies that there are continuum many distinct
relations up to Borel bireducibility.

This last result holds also in the context of $\leq_c$ and $\leq_e$
reducibility, since the arguments that Adams and Kechris use to
establish non-reducibility are measure-theoretic.  We presently define
a class of infinite time computable relations which we propose is the
correct analogue of the countable Borel equivalence relations, and
investigate the corresponding generalizations of the remaining
results.  The idea comes from the classical proof of the maximality of
$E_\infty$, which hinges on the following characterization of the
countable Borel equivalence relations.  Namely, $E$ is a countable
Borel equivalence relation if and only if it admits a \emph{Borel
  enumeration}, that is, a Borel function $f$ such that $f(x)$ codes
an enumeration of $[x]_E$, for all $x$.  (This characterization is an
immediate consequence of the Lusin-Novikov theorem from descriptive
set theory.)  Generalizing this, we say that the equivalence relation
$E$ is (infinite time) \emph{enumerable} if there exists an infinite
time computable function $f$ such that $f(x)$ codes an enumeration of
$[x]_E$, for all $x$.  The \emph{eventually enumerable} equivalence
relations are defined analogously.  This is a worthwhile
generalization; for instance the relation defined by $x\equiv_{\sf
  hyp}y$ iff $x$ and $y$ are hyperarithmetic in one another is
enumerable but not Borel.

Since we have said that the maximality of $E_\infty$ depends on the
above characterization of the countable Borel equivalence relations,
and since we have defined the enumerable and eventually enumerable
equivalence relations in the analogous way, the proof of maximality of
$E_\infty$ in the Borel context yields the same in our context.

\begin{theorem}
  \label{Theorem.Einfty}
  $E_\infty$ is $\leq_c$-greatest among the enumerable relations, and
  $\leq_e$-greatest among the eventually enumerable relations.
\end{theorem}

Perhaps surprisingly, one can establish the minimality of $=$ as well.

\begin{theorem}
  \label{Theorem.Silver}
  $=$ is reducible to every eventually enumerable equivalence relation
  by a continuous function.
\end{theorem}

This result is an immediate consequence of the fact (due originally to
Welch) that there exists a perfect set of $\oequiv_{e\infty}$-classes.
(Here, $\oequiv_{e\infty}$ denotes the \emph{eventual degree}
relation, which is defined analogously to $\oequiv_\infty$.)  The idea
of Welch's proof is to use the theory of forcing over $L_\Sigma$ to
obtain a perfect set of mutually generic Cohen reals, and then argue
that this set does the job.  To see that Theorem~\ref{Theorem.Silver}
follows, observe that every eventually enumerable relation $E$ is
contained (as a set of pairs) in the relation $\oequiv_{e\infty}$.
Hence there exists a perfect set of $E$-classes, and it follows that
there is a continuous reduction from $=$ to $E$.

Finally, we have been unable to establish the minimality of $E_0$ over
the equality relation, and we leave this as a question.  It is hoped
that methods similar to the proof of Theorem~\ref{Theorem.Silver} will
provide an answer.

\begin{question}
  Is it true of every enumerable equivalence relation $E$ that either
  $E$ is reducible to $=$ or else $E_0$ is reducible to $E$?
\end{question}

\bibliographystyle{alpha}
\bibliography{MathBiblio,HamkinsBiblio,ittm}

\newcommand{\etalchar}[1]{$^{#1}$}
\begin{thebibliography}{HMSW07}

\bibitem[CFK{\etalchar{+}}10]{koepke-itrm}
Merlin Carl, Tim Fischbach, Peter Koepke, Russell Miller, Miriam Nasfi, and
  Gregor Weckbecker.
\newblock The basic theory of infinite time register machines.
\newblock {\em Arch. Math. Logic}, 49(2):249--273, 2010.

\bibitem[CH11]{CoskeyHamkins:InfiniteTimeComputableEquivalenceRelations}
Sam Coskey and Joel~David Hamkins.
\newblock Infinite time computable equivalence relations.
\newblock {\em Notre Dame Journal of Formal Logic}, 2011.
\newblock to appear.

\bibitem[DHS05]{DeolalikarHamkinsSchindler2005:NPcoNP}
Vinay Deolalikar, Joel~David Hamkins, and Ralf-Dieter Schindler.
\newblock {${\rm P}\not={\rm NP}\cap {\rm co-NP}$}\/ for infinite time turing
  machines.
\newblock {\em Journal of Logic and Computation}, 15(5):577--592, October 2005.

\bibitem[FS89]{friedman}
Harvey Friedman and Lee Stanley.
\newblock A {B}orel reducibility theory for classes of countable structures.
\newblock {\em J. Symbolic Logic}, 54(3):894--914, 1989.

\bibitem[Ham02]{Hamkins2002:Turing}
Joel~David Hamkins.
\newblock Infinite time turing machines.
\newblock {\em Minds and Machines}, 12(4):521--539, 2002.
\newblock (special issue devoted to hypercomputation).

\bibitem[Ham04]{Hamkins2004:SupertaskComputation}
Joel~David Hamkins.
\newblock Supertask computation.
\newblock In Boris~Piwinger Benedikt~L{\"o}we and Thoralf R{\"a}sch, editors,
  {\em Classical and New Paradigms of Computation and their Complexity
  Hierarchies}, volume~23 of {\em Trends in Logic}, pages 141--158. Kluwer
  Academic Publishers, 2004.
\newblock Papers of the conference ``Foundations of the Formal Sciences III''
  held in Vienna, September 21-24, 2001.

\bibitem[Ham05]{Hamkins2005:InfinitaryComputabilityWithITTM}
Joel~David Hamkins.
\newblock Infinitary computability with infinite time {Turing} machines.
\newblock In Barry~S. Cooper and Benedikt {L\"owe}, editors, {\em New
  Computational Paradigms}, volume 3526 of {\em LNCS}, Amsterdam, June 8-12
  2005. CiE, Springer-Verlag.

\bibitem[Ham07]{Hamkins2007:ASurveyOfInfiniteTimeTuringMachines}
Joel~David Hamkins.
\newblock A survey of infinite time {Turing} machines.
\newblock In {J\'er\^ ome} Durand-Lose and Maurice Margenstern, editors, {\em
  Machines, Computations, and Universality - 5th International Conference MCU
  2007}, volume 4664 of {\em Lecture Notes in Computer Science}, Orleans,
  France, 2007.

\bibitem[Hjo00]{hjorthbook}
Greg Hjorth.
\newblock {\em Classification and orbit equivalence relations}, volume~75 of
  {\em Mathematical Surveys and Monographs}.
\newblock American Mathematical Society, Providence, RI, 2000.

\bibitem[HK96]{hjorthkechris}
Greg Hjorth and Alexander~S. Kechris.
\newblock Borel equivalence relations and classifications of countable models.
\newblock {\em Ann. Pure Appl. Logic}, 82(3):221--272, 1996.

\bibitem[HL00]{HamkinsLewis2000:InfiniteTimeTM}
Joel~David Hamkins and Andy Lewis.
\newblock Infinite time {T}uring machines.
\newblock {\em J. Symbolic Logic}, 65(2):567--604, 2000.

\bibitem[HL02]{HamkinsLewis2002:PostProblem}
Joel~David Hamkins and Andy Lewis.
\newblock Post's problem for supertasks has both positive and negative
  solutions.
\newblock {\em Archive for Mathematical Logic}, 41(6):507--523, 2002.

\bibitem[HLM07]{HamkinsLinetskyMiller2007:ComplexityOfQuicklyDecidableORMSets}
Joel~David Hamkins, David Linetsky, and Russell Miller.
\newblock The complexity of quickly decidable {ORM}-decidable sets.
\newblock In Barry Cooper, Benedikt {L\"owe}, and Andrea Sorbi, editors, {\em
  Computation and Logic in the Real World - Third Conference of Computability
  in Europe CiE 2007}, volume 4497 of {\em Proceedings, Lecture Notes in
  Computer Science}, Siena, Italy, 2007.

\bibitem[HM07]{HamkinsMiller2007:PostsProblemForORMs}
Joel~David Hamkins and Russell Miller.
\newblock Post's problem for ordinal register machines.
\newblock In Barry Cooper, Benedikt {L\"owe}, and Andrea Sorbi, editors, {\em
  Computation and Logic in the Real World - Third Conference of Computability
  in Europe CiE 2007}, volume 4497 of {\em Proceedings, Lecture Notes in
  Computer Science}, Siena, Italy, 2007.

\bibitem[HM09]{HamkinsMiller:PostsProblemForORMsExplicitApproach}
Joel~David Hamkins and Russell Miller.
\newblock Post's problem for ordinal register machines: an explicit approach.
\newblock {\em Annals of Pure and Applied Logic}, 160(3):302--309, 2009.

\bibitem[HMSW07]{HamkinsMillerSeaboldWarner2007:InfiniteTimeComputableModelThe%
ory}
J.~D. Hamkins, R.~Miller, D.~Seabold, and S.~Warner.
\newblock Infinite time computable model theory.
\newblock In S.B.\ Cooper, Benedikt L\"owe, and Andrea Sorbi, editors, {\em New
  Computational Paradigms: Changing Conceptions of What is Computable}, pages
  521--557. Springer, 2007.

\bibitem[HS01]{HamkinsSeabold2001:OneTape}
Joel~David Hamkins and Daniel Seabold.
\newblock Infinite time {Turing} machines with only one tape.
\newblock {\em Mathematical Logic Quarterly}, 47(2):271--287, 2001.

\bibitem[HW03]{HamkinsWelch2003:Pf=NPf}
Joel~David Hamkins and Philip Welch.
\newblock {$P^f\not=NP^f$} for almost all $f$.
\newblock {\em Mathematical Logic Quarterly}, 49(5):536--540, 2003.

\bibitem[KK06]{koepke-ordinal}
Peter Koepke and Martin Koerwien.
\newblock Ordinal computations.
\newblock {\em Math. Structures Comput. Sci.}, 16(5):867--884, 2006.

\bibitem[Koe05]{Koepke2005:TuringComputationsOnOrdinals}
Peter Koepke.
\newblock Turing computations on ordinals.
\newblock {\em Bulletin of Symbolic Logic}, 11(3):377--397, September 2005.

\bibitem[KS06]{Koepke-Siders}
Peter Koepke and Ryan Siders.
\newblock {Register computations on ordinals}.
\newblock {\em submitted to: Archive for Mathematical Logic}, 2006.

\bibitem[KS09]{koepke-recursion}
Peter Koepke and Benjamin Seyfferth.
\newblock Ordinal machines and admissible recursion theory.
\newblock {\em Ann. Pure Appl. Logic}, 160(3):310--318, 2009.

\bibitem[L\"01]{Loewe2001:RevisionSequences}
Benedikt L\"owe.
\newblock Revision sequences and computers with an infinite amount of time.
\newblock {\em Logic Comput.}, 11(1):25--40, 2001.

\bibitem[LM04]{LenziMonteleone2004:OnFixpointArithmeticAndITTM}
Giacomo Lenzi and Erich Monteleone.
\newblock On fixpoint arithmetic and infinite time turing machines.
\newblock {\em Information Processing Letters}, 91(3):121--128, 2004.

\bibitem[Sch03]{Schindler:Pnot=NP}
Ralf-Dieter Schindler.
\newblock {P $\neq$ NP} for infinite time {Turing} machines.
\newblock {\em Monatshefte f\"{u}r Mathematik}, 139(4):335--340, 2003.

\bibitem[ST11]{thomas_workshop}
Scott Schneider and Simon Thomas.
\newblock Countable {B}orel equivalence relations.
\newblock In {\em Appalachian Set Theory 2006--2010}, 2011.

\bibitem[Wel]{Welch:OnAQuestionOfDeolalikarHamkinsSchindler}
Philip Welch.
\newblock On a question of {Deolalikar}, {Hamkins} and {Schindler}.
\newblock available on the author's web page at
  http://www2.maths.bris.ac.uk/$\sim$mapdw/dhs.ps.

\bibitem[Wel99]{Welch99:FriedmansTrick}
Philip Welch.
\newblock Friedman's trick: Minimality arguments in the infinite time {Turing}
  degrees.
\newblock {\em in ``Sets and Proofs'', Proceedings ASL Logic Colloquium},
  258:425--436, 1999.

\bibitem[Wel00a]{Welch2000:Eventually}
Philip Welch.
\newblock Eventually infinite time {Turing} machine degrees: Infinite time
  decidable reals.
\newblock {\em Journal of Symbolic Logic}, 65(3):1193--1203, 2000.

\bibitem[Wel00b]{Welch2000:LengthsOfITTM}
Philip Welch.
\newblock The lengths of infinite time {Turing} machine computations.
\newblock {\em Bulletin of the London Mathematical Society}, 32(2):129--136,
  2000.

\bibitem[Wel05]{Welch2005:ActionOfOneTapeMachines}
Philip Welch.
\newblock The transfinite action of 1 tape {Turing} machines.
\newblock In Barry~S. Cooper and Benedikt {L\"owe}, editors, {\em New
  Computational Paradigms}, volume 3526 of {\em LNCS}, Amsterdam, June 8-12
  2005. CiE, Springer-Verlag.

\end{thebibliography}

\end{document}